\documentclass[12pt,a4paper]{article}

\usepackage{amsthm,amsfonts,amsmath,amssymb,bbm,times1,color}
\usepackage{mysec,euscript,mathrsfs}
\usepackage{hyperref}
\usepackage{lineno}

\usepackage{cancel}
\usepackage[normalem]{ulem}

\newcommand{\PP}{\mathop{{\mathbb{P}}{}}\nolimits}
\newcommand{\EE}{\mathop{{\mathbb{E}}{}\myp}\nolimits}
\newcommand{\const}{\mathrm{const}}
\newcommand{\rd}{{\rm d}}
\newcommand{\re}{{\rm e}}

\newcommand{\g}{\mathrm{g}}

\newcommand{\bfOne}{\mathbbm{1}}

\newcommand{\myp}{\mbox{$\:\!$}}
\newcommand{\mypp}{\mbox{$\;\!$}}
\newcommand{\myn}{\mbox{$\:\!\!$}}
\newcommand{\mynn}{\mbox{$\;\!\!$}}

\newcommand{\RR}{{\mathbb R}}

\newcommand{\NN}{{\mathbb N}}
\newcommand{\ZZ}{{\mathbb Z}}
\newcommand{\QQ}{{\mathbb Q}}
\newcommand{\T}{T\kern-.05pc{}}

\newcommand{\calA}{{\EuScript A}}
\newcommand{\calB}{{\EuScript B}}
\newcommand{\calD}{{\EuScript D}}
\newcommand{\calF}{{\EuScript F}}
\newcommand{\calG}{{\EuScript G}}

\newcommand{\calI}{{\EuScript I}}

\numberwithin{equation}{section}

\newtheorem{theorem}{Theorem}[section]
\newtheorem{lemma}[theorem]{Lemma}
\newtheorem{proposition}[theorem]{Proposition}

\theoremstyle{remark}
\newtheorem{remark}{Remark}[section]

\theoremstyle{definition}
\newtheorem{example}{Example}[section]
\newtheorem{assumption}{Assumption}[section]
\newtheorem{definition}{Definition}[section]

\makeatletter
\def\MR#1{\href{http://www.ams.org/mathscinet-getitem?mr=#1}{MR#1}}

\makeatother

\begin{document}

\title{On bounded continuous solutions of the archetypal equation with rescaling}
\author{Leonid V.~Bogachev\myp$^{\rm a}$, Gregory Derfel\myp$^{\rm b}$ and Stanislav A.~Molchanov\myp$^{\rm c}$}
\date{\small $^{\,\rm a}$\myp Department of Statistics,
School of Mathematics, University of
Leeds,\\ Leeds LS2 9JT, UK. E-mail: {\tt L.V.Bogachev@leeds.ac.uk}\\[.2pc]
$^{\,\rm b}$\myp Department of Mathematics, Ben Gurion University of
the Negev,\\
Be'er Sheva 84105, Israel. E-mail: {\tt derfel@math.bgu.ac.il}\\[.2pc]
$^{\,\rm c}$\myp Department of Mathematics, University of North
Carolina at Charlotte,\\ Charlotte NC\,28223, USA. E-mail: {\tt
smolchan@uncc.edu}}


\maketitle

\begin{abstract}

\medskip\noindent
The `archetypal' equation with rescaling is given by
$y(x)=\iint_{\RR^2} y(a(x-b))\,\mu(\rd{a},\rd{b})$ ($x\in\RR$),
where $\mu$ is a probability measure; equivalently,
$y(x)=\EE\{y(\alpha\myp(x-\beta))\}$, with random $\alpha,\beta$ and
$\EE$ denoting expectation. Examples include: (i) functional
equation $y(x)=\sum_{i} p_{i}\mypp y(a_i(x-b_i))$; (ii)
functional-differential (`pantograph') equation $y'(x)+y(x)
=\sum_{i} p_{i}\mypp y(a_i(x-c_i))$ ($p_{i}>0$, $\sum_{i} p_{i}=1$).
Interpreting solutions $y(x)$ as harmonic functions of the
associated Markov chain $(X_n)$, we obtain Liouville-type results
asserting that any bounded continuous solution is constant.
In particular, in the `critical' case $\EE\{\ln\myn|\alpha|\}=0$
such a theorem holds subject to uniform continuity of $y(x)$; the
latter is guaranteed under mild regularity assumptions on
$\beta$, satisfied e.g.\ for the pantograph equation (ii).
For equation (i) with $a_i=q^{m_i}$ \myp(
$m_i\in\ZZ$, $\sum_i p_i\myp m_i=0$), the result can be proved
without the uniform continuity assumption.
The proofs utilize the iterated equation
$y(x)=\EE\{y(X_\tau)\myp|\myp X_0=x\}$ (with a suitable stopping
time $\tau$) due to Doob's optional stopping theorem applied to the
martingale~$y(X_n)$.

\vspace{.4pc}\noindent \emph{Keywords}\/: Functional \&
functional-differential equations; pantograph equation; Markov
chain; harmonic function; martingale; stopping time

\vspace{.4pc}\noindent \emph{2010 MSC}:
Primary: 39B05; Secondary: 34K06, 39A22, 60G42, 60J05. \hfill{}

\end{abstract}


\section{Introduction}\label{sec:intro}
\subsection{The archetypal equation}\label{sec:intro1}
This paper concerns the equation with rescaling (referred to as
`archetypal') of the form
\begin{equation}\label{eq:eq0}
y(x)=\iint_{\RR^2} y(a(x-b))\,\mu(\rd{a},\rd{b}),\qquad x\in\RR,
\end{equation}
where $\mu(\rd{a},\rd{b})$ is a probability measure on $\RR^2$. The
integral in \eqref{eq:eq0} has the meaning of expectation with
respect to a random pair $(\alpha,\beta)$ with distribution
$\PP\{(\alpha,\beta)\in \rd{a}\times\rd{b}\}=\mu(\rd{a},\rd{b})$;
thus, equation \eqref{eq:eq0} can be written in compact form as
\begin{equation}\label{eq:arch}
y(x)=\EE\{y(\alpha\myp(x-\beta))\},\qquad x\in\RR.
\end{equation}

As has been observed by Derfel \cite{Der} (and will be illustrated
below in \S\ref{sec:subclasses}), this equation is a rich source of
various functional and functional-differential equations with
rescaling, specified by a suitable choice of the probability measure
$\mu$ (i.e.\ the distribution of $(\alpha,\beta)$). It is for that
reason that we propose to call \eqref{eq:arch} (as well as its
integral counterpart~\eqref{eq:eq0}) the \emph{archetypal
equation}~(\emph{AE}). The study of this equation allows one to
enhance and unify earlier results for particular subclasses of
equations with rescaling, while making the analysis more transparent
and efficient (cf.\ \cite{BDMO}).

Noting that any function $y(x)\equiv\const$ satisfies the AE
\eqref{eq:arch}, it is natural to investigate whether there are any
\textit{non-trivial} (i.e.\ non-constant) bounded continuous (b.c.)\
solutions. Such a question naturally arises in the context of
functional and functional-differential equations with rescaling,
where the possible existence of bounded solutions (e.g.\ periodic,
almost periodic, compactly supported, etc.) is of major interest in
physical and other applications (see e.g.\
\cite{Cavaretta,Rvachev,Schilling,Spiridonov}). Solutions under
study may also be bounded by nature, e.g.\ representing a ruin
probability as a function of initial capital \cite{Ferguson,Mahler}.
On the other hand, confining oneself to bounded solutions may be
considered as a first step towards a full description of the
asymptotic behaviour of solutions.

Thus, the goal of the present paper is to give conditions on the
distribution $\mu$ of the random coefficients $(\alpha,\beta)$,
under which any b.c.-solution of equation \eqref{eq:arch} is
constant on $\RR$. For shorthand, we refer to statements of this
kind as \emph{Liouville-type} theorems by analogy with similar
results in complex analysis and harmonic function theory,
bearing in mind that $y(x)$ in \eqref{eq:arch} is a weighted average
of other values, thus resembling the usual harmonic function. More
details highlighting the pertinence of `harmonicity' in the context
of equation \eqref{eq:arch} are provided in~\S\ref{sec:1.4}.

\begin{remark}
Continuity of $y(x)$ (or some other regularity assumption) is needed
to avoid pathological solutions, as is well known in the theory of
functional equations (cf.\ \cite[Ch.\,2]{Aczel}). For example, all
b.c.-solutions of the equation $y(x)=p\mypp
y\bigl(\frac12(x+1)\bigr)+(1-p)\myp y\bigl(\frac12(x-1)\bigr)$
($0\le p\le1$) are constant by Theorem \ref{th:K<>0}(a) stated
below, but if the continuity requirement is dropped then one can
easily construct other bounded solutions, e.g.\ the Dirichlet
function $y(x)=\bfOne_{\QQ}(x)$ (i.e.\ $y(x)=1$ if $x$ is rational
and $y(x)=0$ otherwise), which is everywhere discontinuous.
\end{remark}

\subsection{Some subclasses of the archetypal equation; historical
remarks}\label{sec:subclasses}

Before outlining our results, we illustrate the remarkable capacity
of
equation
\eqref{eq:arch} justifying the name `archetypal'. General surveys of
functional and functional-differential equations with rescaling are
found in Derfel \cite{Der1} and Baron \& Jarczyk \cite{Baron}, both
with extensive bibliographies.

\subsubsection{Functional equations and self-similar measures.}\label{sec:1.2.1}
To start with, in the simplest case $\alpha\equiv 1$ equations
\eqref{eq:eq0}, \eqref{eq:arch} are reduced to
\begin{equation}\label{eq:CD}
y(x)=\int_\RR y(x-t)\,\mu_\beta(\rd{t})\quad\Longleftrightarrow\quad
y(x)=\EE\{y(x-\beta)\},
\end{equation}
where $\mu_\beta(\rd{t}):=\PP(\beta\in\rd{t})$. This equation
(sometimes called the \emph{integrated Cauchy functional equation}
\cite{Rao}) plays a central role in potential theory and harmonic
analysis on groups \cite{Chu,Meyer}, and is also prominent in
probability theory in relation to renewal theorems
\cite[\S{}XI.9]{Feller2}, Markov chains
\cite[Ch.\;5]{Revuz}, queues
\cite[\S{}III.6]{Asmussen}, characterization of probability
distributions \cite[Ch.\;2]{Rao}, etc. A Liouville-type result in
this case is rendered by the celebrated \emph{Choquet--Deny
theorem}~\cite{CD} (see also
\cite{Rao} and references therein).

Note that equation \eqref{eq:CD} can be written in the
convolution\footnote{The convolution between function $y(x)$ and
measure $\sigma$ in $\RR$ is defined as
$y\star\sigma(x):=\int_{\RR}y(x-t)\,\sigma(\rd t)$.} form
$y=y\star\mu_\beta$. More generally, if $\alpha$ has a discrete
distribution (with atoms $a_i$ and masses $p_i$) then, denoting by
$\mu_\beta^{\myp i}$ the conditional distribution of $\beta$ given
$\alpha=a_i$, the AE
\eqref{eq:arch} is conveniently expressed in convolutions,
\begin{equation*}
y(x)=\sum\nolimits_{i} p_i\mypp y(a_ix)\star\mu_\beta^{\myp i}\mypp.
\end{equation*}
For a purely discrete measure $\mu$, with atoms $(a_i,b_i)$ and
masses $p_i=\PP(\alpha=a_i,\mypp\beta=b_i)$, the AE \eqref{eq:arch}
specializes to
\begin{equation}\label{eq:pure1}
y(x)=\sum\nolimits_{i} p_i\mypp y(a_i(x-b_{i})).
\end{equation}
If all $a_i>1$ then \eqref{eq:pure1} is an example of
\emph{Hutchinson's equation} \cite{Hutchinson} for the distribution
function of a self-similar probability measure which is invariant
under a family of contractions (here, affine transformations
$x\mapsto b_i+x/a_i$). An important subclass of \eqref{eq:pure1},
with $a_i\equiv a>1$, is exemplified by
\begin{equation*}
y(x)=\tfrac12\mypp y(a\myp(x+1))+\tfrac12\mypp y(a\myp(x-1)).
\end{equation*}
This equation describes the (self-similar) distribution function of
the random series\break $\sum_{n=0}^\infty \pm\myp a^{-n}$, where
the signs are chosen independently with probability $\frac12$.
Characterization of this distribution for different $a>1$ is the
topical \emph{Bernoulli convolutions problem}~\cite{Peres}.

Returning to equation \eqref{eq:pure1} with $a_i\equiv a>1$ the
density $z(x):=y'(x)$ (if it exists) satisfies
\begin{equation}\label{eq:two-scale'}
z(x)=a\sum\nolimits_i p_i\mypp z(a\myp(x-b_i)),
\end{equation}
often called the \emph{two-scale difference equation} or
\emph{refinement equation}~\cite{DL}. Construction of compactly
supported continuous solutions of \eqref{eq:two-scale'} plays a
crucial role in wavelet theory \cite{D10,Strang} and also in
subdivision schemes and curve design \cite{Cavaretta,DerDynLev},
which is a rapidly growing branch of approximation theory. A special
version of \eqref{eq:two-scale'} known as \emph{Schilling's
equation}
\begin{equation*}
z(x)=\alpha\mypp\bigl[\tfrac14\myp z(\alpha\myp x+1) + \tfrac12\myp
z(\alpha\myp x) + \tfrac14\myp z(\alpha\myp x-1)\bigr]
\end{equation*}
arises in solid state physics in relation to spatially chaotic
structures in amorphous materials
\cite[
p.\;230]{Schilling}, where the existence of compactly supported
continuous solutions is again of major interest; see \cite{DerSchi}
for a full characterization of this problem in terms of arithmetical
properties of~$\alpha$.

\subsubsection{Functional-differential equations.}\label{sec:pantograph}

Let us now turn to the situation where the distribution of $\beta$
conditioned on $\alpha$ is absolutely continuous (i.e.\ has a
density). It appears that for certain simple densities the AE
\eqref{eq:arch} produces some well-known functional-differential
equations. An important example is the celebrated \emph{pantograph
equation}, introduced by Ockendon \& Tayler \cite{OT} as a
mathematical model of the overhead current collection system on an
electric locomotive.\footnote{The term `pantograph equation' was
coined by Iserles~\cite{Iserles}.} In its classical
(one-dimensional) form the pantograph equation reads
\begin{equation}\label{eq:pant}
y'(x)=c_0\myp y(x)+c_1 y(\alpha\myp x).
\end{equation}
This equation and its ramifications have emerged in a striking range
of applications, including number theory
\cite{Mahler}, astrophysics
\cite{Chandra},
queues \& risk theory
\cite{Gaver},
stochastic games \cite{Ferguson},
quantum theory \cite{Spiridonov},
and population dynamics
\cite{HallWake1},
The common feature of all such examples is some self-similarity of
the system under study. Thorough asymptotic analysis of equation
\eqref{eq:pant} was given by Kato \& McLeod \cite{KM}.
A more general first-order pantograph equation (with matrix
coefficients, and also allowing for a term with a rescaled
derivative) was studied by Iserles \cite{Iserles}, where a fine
geometric structure of almost-periodic solutions was also described.
Further developments include analysis in the complex domain
\cite{DI}, higher-order equations \cite{vanBrunt,Iserles}, and
stochastic versions~\cite{Fan}. Among recent important analytic
results is a proof by da Costa \emph{et al.}\ \cite{DaCosta} of the
unimodality of solutions which plays a significant role in medical
imaging of tumours~\cite{Sueb}.

A \emph{balanced} version of the pantograph equation is given by
(see \cite{BDMO,Der})
\begin{equation}\label{eq:pant-gen}
y'(x)+y(x)=\sum\nolimits_{i} p_i\mypp y(a_i\myp (x-c_i)),\qquad
p_i>0,\quad \sum\nolimits_{i} p_i=1.
\end{equation}
As explained in
\S\myp\ref{sec:A1}, equation \eqref{eq:pant-gen} is essentially
equivalent to the AE \eqref{eq:arch} whereby $\alpha$ is discrete,
with $\PP(\alpha=a_i)=p_i$, and $\beta$ conditioned on $\alpha=a_i$
has the unit exponential distribution on $(c_i,\infty)$, with the
density function $\re^{c_i-t}\bfOne_{(c_i,\infty)}(t)$. The
discreteness of $\alpha$ is not significant here, and a similar
connection with the AE holds for more general integro-differential
equations (cf.\ \cite{Iserles-Liu})
\begin{equation}\label{eq:pant-gen-alpha}
y'(x)+y(x)=\EE\{y(\alpha\myp (x-\gamma)\}\equiv\iint_{\RR^2}
y(a\myp(x-c))\,\mu_{\alpha,\gamma}(\rd{a},\rd{c}),
\end{equation}
where $\gamma$ is a random variable and
$\mu_{\alpha,\gamma}(\rd{a},\rd{c})=\PP(\alpha\in\rd{a},\gamma\in\rd{c})$
is the distribution of $(\alpha,\gamma)$. Higher-order pantograph
equations can also be deduced from the AE,
e.g.
$$
-y''(x)+y(x)= \sum\nolimits_{i} p_i\mypp y(a_i\myp (x-c_i)),\qquad
p_i>0,\quad \sum\nolimits_{i} p_i=1,
$$
and more generally (cf.\ \cite{BDMO})
$$
C_2\mypp y''(x)+C_1\myp y'(x)+y(x)=\EE\{y(\alpha\myp
(x-\gamma)\}\qquad (C_1,C_2\in\RR,\ \ C_1^2-4\myp C_2 \ge 0).
$$

For an example of a different kind, take $\alpha\equiv 2$ and assume
that $\beta$ has the unform distribution on $[-\frac12,\frac12]$,
then equation \eqref{eq:arch} is reduced to
\begin{equation}\label{eq:Rvachev1}
y(x)=\int_{x-1/2}^{x+1/2} y(2\myp u)\,\rd{u}.
\end{equation}
Differentiating \eqref{eq:Rvachev1}, for $z(x):=y'(x)$ we obtain
\emph{Rvachev's equation} \cite{Rvachev}
\begin{equation}\label{eq:Rvachev}
z'(x)=2\bigl[z(2x+1)-z(2x-1)\bigr].
\end{equation}
A compactly supported solution of \eqref{eq:Rvachev} (called the
`\emph{up}-function') and its generalizations (unified under the
name \emph{atomic functions}) have extensive applications in
approximation theory (see \cite{DerDynLev,Rvachev} and references
therein); all such functions can be obtained as solutions of
suitable versions of the AE \eqref{eq:arch}
(see~\cite{Der}).

\subsection{Main results}\label{sec:intro2}

Let us summarize our results. First, certain degenerate cases
warrant a separate analysis but need to be excluded in general
theory, namely: (i) $\alpha=0$ with positive probability; (ii)
$|\alpha|\equiv 1$; and (iii) $\alpha\myp(c-\beta)\equiv c$ for some
$c\in\RR$ (\emph{resonance}). Note that (ii) includes the case
$\alpha\equiv1$ settled in the Choquet--Deny theorem mentioned in
\S\ref{sec:1.2.1}; in \S\myp\ref{sec:CD} we generalize this result
(Theorem~\ref{th:CD}). As for cases (i) and (iii), a Liouville
theorem holds here unconditionally, which is easy to prove for (i),
analytically and probabilistically alike (see
Theorem~\ref{th:alpha=0}). In the resonance case (iii), the proof is
more involved relying heavily on the Choquet--Deny theorem (see
\S\myp\ref{sec:resonance}), but the result itself is quite lucid and
appealing.

In the non-degenerate situation, existence of non-trivial
b.c.-solutions is essentially governed by the sign of
$K:=\iint_{\RR^2}\ln\myn|a|\,\mu(\rd{}a\times\rd{}b)=\EE\{\ln\myn|\alpha|\}$.
More precisely, one can prove (see~\cite{BDM}) the following
dichotomy between the \emph{subcritical} ($K<0$) and
\emph{supercritical} ($K>0$) regimes.

\begin{theorem}\label{th:K<>0} Suppose that $K=\EE\{\ln \mynn|\alpha|\}$
is finite and $\EE\{\ln\max(|\beta|,1)\}<\infty$.
\begin{itemize}
\item[\rm (a)] If $K<0$ then any
b.c.-solution of the AE\/ \eqref{eq:arch} is constant.
\item[\rm (b)]
If $K>0$ and $\alpha>0$, then there is a b.c.-solution of\/
\eqref{eq:arch} given by the distribution function
$F_\varUpsilon(x):=\PP(\varUpsilon\le x)$, where $\varUpsilon:
=\sum_{n=1}^\infty \beta_n\prod_{j=1}^{n-1}\alpha_j^{-1}$ and
$\{(\alpha_n,\beta_n)\}$ is a sequence of independent identically
distributed (i.i.d.)\ random pairs with the same distribution $\mu$
as $(\alpha,\beta)$.
\end{itemize}
\end{theorem}

\begin{remark}
Almost sure (a.s.)\ convergence of the random series $\varUpsilon$
(for $\alpha\ne0$) and continuity of $F_\varUpsilon(x)$ on $\RR$
were proved by Grintsevichyus \cite{Gr}.
\end{remark}
\begin{remark}
The result of Theorem~\ref{th:K<>0} was obtained by Derfel
\cite{Der}
under a stronger moment condition $\EE\{|\beta|\}<\infty$ and only
for $\alpha>0$ (which is essential in~(b) but not in~(a)); however,
his arguments hold in the general case with minor changes.
\end{remark}

\begin{remark}
In contrast with the subcritical case $K<0$, which is insensitive to
the sign of $\alpha$ (see Theorem~\ref{th:K<>0}(a)), the
supercritical case $K>0$ is more delicate: if $\PP(\alpha<0)>0$ then
$y=F_\varUpsilon(x)$ is no longer a solution of the
AE~\eqref{eq:arch}; e.g.\ if $\alpha<0$ (a.s.)\ then this function
satisfies the equation $y(x)=1-\EE\{y(\alpha\myp(x-\beta))\}$ (cf.\
\cite[Eq.\,(5)]{Gr}). Moreover, one can prove \cite{BDM} that any
bounded solution of \eqref{eq:arch} with limits at $\pm\infty$ is
constant; thus, any non-trivial solution must be oscillating, which
is drastically different from the case $\alpha>0$~(a.s.).
\end{remark}

The \emph{critical case} $K=0$ is much more challenging, and it has
remained open since~\cite{Der}. More recently, for a pantograph
equation \eqref{eq:pant-gen-alpha} without shift (i.e.\
$\gamma\equiv0$) and some second-order extensions, a Liouville
theorem in the case $K=0$ was established by Bogachev \emph{et
al.}~\cite{BDMO}. In the present paper, we prove the following
general result (cf.\ Theorem \ref{th:K=0} below).

\begin{theorem}\label{th:K=0_intro}
Assume that\/ $\PP(|\alpha|\ne 1\}>0$ and\/
$K=\EE\{\ln\myn|\alpha|\}=0$. If\/ $y(x)$ is a bounded solution of
the AE\/ \eqref{eq:arch} which is uniformly continuous on $\RR$,
then it is constant.
\end{theorem}

Although an unwanted restriction, the uniform continuity assumption
can be shown to be satisfied provided there exists the probability
density of $\beta$ conditioned on $\alpha$
(Theorem~\ref{th:uniform}). An alternative criterion tailored to the
model $\beta=\gamma+\xi$ with $\xi$ independent of $(\alpha,\gamma)$
(Theorem~\ref{th:uniform2}) is applicable to a large class of
examples including the pantograph equation \eqref{eq:pant-gen-alpha}
and its generalizations (\S\myp\ref{sec:A1}). As a consequence, we
obtain a Liouville theorem for the general (balanced) pantograph
equation in the critical case (cf.\ Theorem~\ref{th:pantograph}),
significantly extending the result of Bogachev \emph{et
al.}~\cite{BDMO}. In particular, for the first-order pantograph
equation \eqref{eq:pant-gen-alpha} we have

\begin{theorem}\label{th:pantograph1}
If\/ $\PP(|\alpha|\ne 1\}>0$ and $K=\EE\{\ln\myn|\alpha|\}=0$, then
any bounded solution of \eqref{eq:pant-gen-alpha} is constant.
\end{theorem}

\begin{remark} As explained in \S\myp\ref{sec:3.3} below,
Theorem \ref{th:pantograph1} extends to the case $|\alpha|=1$
(a.s.)\ by virtue of the generalized Choquet--Deny theorem proved in
\S\myp\ref{sec:CD}.
\end{remark}

On the other hand, for a subclass of functional equations
\eqref{eq:pure1} with multiplicatively commensurable coefficients
$\{a_i\}$ (i.e. $a_i=q^{m_i}$, $q>1$, \,$m_i\in\ZZ$), where the
random shift $\beta$ is discrete and therefore the criteria of
Theorems \ref{th:uniform},~\ref{th:uniform2} do not apply, the
Liouville theorem in the critical case $K=0$ (i.e. $\sum_i p_i\myp
m_i=0$) can be proved by a different method that circumvents the
hypothesis of uniform continuity (Theorem~\ref{th:3}).

\subsection{The method: associated Markov chain and iterations}\label{sec:1.4}

In this subsection, we describe the probabilistic approach to the AE
based on Markov chains and martingales, and introduce some basic
notation and definitions.

Consider a Markov chain $(X_n)$ on $\RR$ defined recursively by
\begin{equation}\label{eq:Xn}
X_n=\alpha_n \myp(X_{n-1}-\beta_n) \quad (n\in\NN), \qquad
X_0=x\in\RR,
\end{equation}
where $\{(\alpha_n,\beta_n)\}_{n\in\NN}$ is a sequence of i.i.d.\
random pairs with the same distribution as $(\alpha,\beta)$ (see
\eqref{eq:arch}). Note that the AE \eqref{eq:arch} is then expressed
as
\begin{equation}\label{eq:Markov,n=1}
y(x)=\EE_x\{y(X_1)\},\qquad x\in\RR,
\end{equation}
where index $x$ in the expectation refers to the initial condition
in \eqref{eq:Xn}. That is to say, any solution of the AE
\eqref{eq:arch} is a \emph{harmonic function} of the Markov chain
$(X_n)$
(cf.\ \cite[
p.\;40]{Revuz}).

\begin{remark}
Bounded harmonic functions play a paramount role in the general
theory of Markov chains (see more details and some references in the
Appendix, \S\myp\ref{sec:A2}).
\end{remark}

\begin{remark}
Stochastic recursion \eqref{eq:Xn} is well known in the literature
as the \emph{random difference equation} (see e.g.\ \cite{Babillot,
Embrechts,Kesten,Vervaat} and further references therein).
\end{remark}

Note that equation \eqref{eq:Markov,n=1} propagates along the Markov
chain $(X_n)$, i.e.\ for any $n\in\NN$
\begin{equation}\label{eq:stop_x_0}
y(x)=\EE_{x}\{y(X_{n})\}, \qquad x\in\RR.
\end{equation}
Equivalently, an integral form of equation \eqref{eq:stop_x_0} can
be obtained by iterating forward  the AE \eqref{eq:eq0}. The
recursion \eqref{eq:Xn} can also be iterated to give explicitly
\begin{equation}\label{eq:Y}
X_{n}=A_n x -D_n\qquad (n\in\NN),
\end{equation}
where
\begin{equation}\label{eq:AY} A_{n}:=\prod_{k=1}^n
\alpha_k, \qquad D_n:=\sum_{k=1}^{n}\beta_{k}\prod_{j=k}^n
\alpha_{j}.
\end{equation}

Recall that $K=\EE\{\ln\myn|\alpha|\}$. In the subcritical case
($K<0$),
for the proof of Theorem~\ref{th:K<>0}(a) it suffices to consider
iterations \eqref{eq:stop_x_0} as $n\to\infty$. Indeed, Kolmogorov's
strong law of large numbers
implies that $S_n:=\sum_{k=1}^n\ln\myn|\alpha_k|\to -\infty$ and
hence $|A_n|=\exp(S_n)\to0$ (a.s.); in view of \eqref{eq:Y} this
indicates that the right-hand side of \eqref{eq:stop_x_0} eventually
becomes $x$-free (see more details in~\cite{BDM}). In the critical
case ($K=0$), the random walk $(S_n)$ is recurrent but none the less
$\liminf_{n\to\infty} S_n=-\infty$ (a.s.); hence, at some random
times $\tau_\epsilon$ we have
$|A_{\tau_\epsilon}|=\exp(S_{\tau_\epsilon})<\epsilon$ (for any
$\epsilon>0$), which can be used to infer that $y(x)\equiv \const$
in a similar fashion as before.

Expanding this idea, our approach to the analysis of equation
\eqref{eq:arch}, first probed in~\cite{BDMO}, is based on replacing
$n$ in the iterated equation \eqref{eq:stop_x_0} by a suitable
\emph{stopping time} $\tau$, defined as a random (integer-valued)
variable such that for any $n\in\NN$ the event $\{\tau\le n\}$ is
determined by $(\alpha_1,\beta_1),\dots,(\alpha_n,\beta_n)$. It
suffices for our purposes to work with `hitting times'
$\tau_B:=\inf\{n\ge1\colon A_n\in B\}\le\infty$,\footnote{Here and
below, we adopt the convention that $\inf\emptyset:=\infty$.} where
$A_n=\alpha_1\cdots\alpha_n$ (see \eqref{eq:AY}) and $B\subset\RR$
is an interval or a single point.

We shall routinely use the following central lemma (where continuity
of $y(x)$ is not required).

\begin{lemma}\label{lm:stop1}
Let $(X_n)$ be the associated Markov chain \eqref{eq:Xn}, and $\tau$
a stopping time such that $\tau<\infty$ a.s. If\/ $y(x)$ is a
bounded solution of the AE \eqref{eq:arch} then it satisfies the
`stopped' equation
\begin{equation}\label{eq:stop_x}
y(x)=\EE_{x}\{y(X_{\tau})\},\qquad x\in\RR.
\end{equation}
\end{lemma}

The crucial fact is that $y(X_n)$ is a
\emph{martingale} (cf.\ \cite[
p.\;43, Proposition~1.8]{Revuz}); indeed, by \eqref{eq:Xn}
\begin{equation}\label{eq:y(X)mart}
\EE_x\myn\{y(X_n)|\myp (\alpha_k,\beta_k),\,k< n\}
=\EE_x\bigl\{y(\alpha_n(X_{n-1}-\beta_n))|\myp X_{n-1}\bigr\}
=y(X_{n-1})\quad \text{(a.s.)},
\end{equation}
which verifies the martingale property
\cite[\S10.3, p.\;94]{Williams}. The lemma then readily follows by
Doob's optional stopping theorem
(e.g.\ \cite[p.\,100, Theorem 10.10(b)]{Williams}). For the sake of
a more self-contained exposition, a direct proof of Lemma
\ref{lm:stop1} is included in the Appendix
(see~\S\myp\ref{sec:lemma}).


\subsubsection*{Layout.}\label{sec:layout} The rest of the paper is
organized as follows. In \S\myp\ref{sec:2} we work out the
degenerate cases mentioned at the beginning of
\S\myp\ref{sec:intro2}, namely: $\PP(\alpha=0)>0$
(\S\myp\ref{sec:alpha=0}), $|\alpha|\equiv 1$ (\S\myp\ref{sec:CD}),
and $\alpha\myp(c-\beta)\equiv c$ (\S\myp\ref{sec:resonance}). In
\S\myp\ref{sec:3.1} we prove our main result for the critical case
(Theorem~\ref{th:K=0}; cf.\ Theorem~\ref{th:K=0_intro}), backed up
in \S\myp\ref{sec:3.2} by simple sufficient conditions for the
uniform continuity of solutions (Theorems \ref{th:uniform}
and~\ref{th:uniform2}). In \S\myp\ref{sec:A1} we explain in detail
the remarkable link between the pantograph equations and the AE,
which enables us to prove a Liouville theorem for the general
(integro-differential) pantograph equation of any order
(Theorem~\ref{th:pantograph}). This is complemented in
\S\myp\ref{sec:3.3} by a Liouville theorem for the functional
equation \eqref{eq:pure1} with $a_i=q^{m_i}$ (Theorem~\ref{th:3}).
The Appendix comprises an elementary proof of Lemma \ref{lm:stop1}
(\S\myp\ref{sec:lemma}) and a brief compendium of basic facts
illuminating the fundamental role of bounded harmonic functions in
the general theory of Markov chains (\S\myp\ref{sec:A2}).

\section{Three degenerate cases}\label{sec:2}
Before embarking on a general discussion of the AE, we need to study
the problem of bounded solutions in certain special cases of
possible values of $\alpha$ and $\beta$, which will be excluded from
consideration thereafter. Recall the notation $A_n:=\prod_{k=1}^n
\alpha_k$ (see~\eqref{eq:AY}).

\subsection{Vanishing of the scaling coefficient}\label{sec:alpha=0}
Let us consider the case where the scaling coefficient $\alpha$ may
take the value zero. Note that continuity of solutions $y(x)$ is not
assumed \textit{a priori}.

\begin{theorem}\label{th:alpha=0}
Suppose $p_0:=\PP(\alpha=0)>0$. Then any bounded solution of the
AE\/ \eqref{eq:arch} is constant on $\RR$.
\end{theorem}

We first give an elementary `analytic' proof of this simple theorem
and then present another proof to illustrate the method based on
Lemma~\ref{lm:stop1}.

\proof[Proof of\/ Theorem \textup{\ref{th:alpha=0}}] With
$y_0(x):=y(x)-y(0)$, equation \eqref{eq:eq0} may be written in the
form
\begin{equation}\label{eq:alpha=0a}
y_0(x)=(1-p_0)\iint_{\RR^2\setminus\{a=0\}}
y_0(a\myp(x-b))\,\tilde{\mu}(\rd{a},\rd{b}),\qquad x\in\RR,
\end{equation}
where $\tilde{\mu}:=(1-p_0)^{-1}\mu$, so that
$\tilde{\mu}(\RR^2\mypp{\setminus}\mypp\{a=0\})=1$. Denoting by
$\|f\|:=\sup_{x\in\RR} |f(x)|$ the sup-norm on $\RR$, from
\eqref{eq:alpha=0a} we obtain
\begin{equation}\label{eq:alpha=0b}
|y_0(x)|\le (1-p_0)\myp\|y_0\| \quad
(x\in\RR)\quad\Longrightarrow\quad \|y_0\|\le (1-p_0)\myp\|y_0\|.
\end{equation}
Since $1-p_0<1$, the second inequality in \eqref{eq:alpha=0b}
immediately implies that $\|y_0\|=0$, and then the first inequality
gives $y_0(x)\equiv 0$, i.e.\ $y(x)\equiv y(0)$, as claimed.
\endproof
\proof[Alternative proof of\/ Theorem \textup{\ref{th:alpha=0}}]
Consider the stopping time $\tau_0:=\inf\{n\ge 1\colon A_n=0\}$.
Note that
\begin{align*}
\PP(\tau_0>n)=\PP(A_n\ne0)&=\PP(\alpha_1\ne0,\dots,\alpha_n\ne0)=(1-p_0)^{n}\to
0\qquad (n\to\infty),
\end{align*}
hence $\tau_0<\infty$ a.s. Now, using the iteration formulas
\eqref{eq:Y}, \eqref{eq:AY}, and noting that $A_{\tau_0}=0$ a.s., by
Lemma \ref{lm:stop1} we obtain
\begin{equation}\label{eq:u=0}
y(x)=\EE_{x}\{y(X_{\tau_0})\}= \EE\{y(x
A_{\tau_0}-D_{\tau_0})\}=\EE\{y(-D_{\tau_0})\},\qquad x\in\RR.
\end{equation}
Since the right-hand side of \eqref{eq:u=0} does not depend on $x$,
it follows that $y(x)=\const$.
\endproof

\subsection{An extension of the Choquet--Deny theorem}\label{sec:CD}

\subsubsection{The classical case $\alpha\equiv1$.} As mentioned in \S\ref{sec:subclasses},
the AE \eqref{eq:arch} with $\alpha\equiv 1$ is reduced to
\begin{equation}\label{eq:CD1}
y(x)=\EE\{y(x-\beta)\},\qquad x\in\RR.
\end{equation}
The famous Choquet--Deny theorem \cite{CD} (cf.\
\cite[
p.\;382, Corollary]{Feller2} or \cite[p.\,161, Theorem~1.3]{Revuz})
asserts that any b.c.-solution of \eqref{eq:CD1} is constant
provided that (the distribution of) the shift $\beta$ is
non-arithmetic, i.e.\ not supported on any set $\lambda\ZZ=\{\lambda
k,\,k\in\ZZ\}$ (with \emph{span} $\lambda\in\RR$).

\begin{remark}\label{rm:UC}
In connection with the uniform continuity condition in Theorem
\ref{th:K=0_intro},
it may be of interest to note that some proofs of the Choquet--Deny
theorem (e.g.\ \cite[p.\;382]{Feller2}) deploy the convolution
$\tilde{y}(x)=y\star \varphi_{0,\myp\sigma^2}(x)$ of a b.c.-solution
$y(x)$ with the density function $\varphi_{0,\myp\sigma^2}(x)$ of
the normal distribution with zero mean and variance $\sigma^2$,
whereby the function $\tilde{y}(x)$ is uniformly continuous and
still satisfies equation \eqref{eq:CD1}. Once it has been proved
that $\tilde{y}(x)$ is constant, this is extended to the original
solution $y(x)$ by taking the limit as $\sigma\to0$.
\end{remark}

A \emph{discrete} version of the Choquet--Deny theorem \cite{CD}
(see also \cite[\S{}XIII.11]{Feller1} or \cite[p.\;276,
Theorem~T1]{Spitzer}) refers to the case where equation
\eqref{eq:CD1} is considered on $\ZZ$ and $\beta$ is integer-valued.
Namely, assume that the smallest additive group containing the set
$\{x\in\ZZ\colon \PP(\beta=x)>0\}$ coincides with $\ZZ$; then the
theorem asserts that $y(x)\equiv\const$ for all $x\in\ZZ$.

In the context of equation \eqref{eq:CD1} on the whole line, this
enables one to give a
full description of b.c.-solutions in the arithmetic case (excluding
the degenerate case $\beta\equiv0$). The next result is essentially
well known in folklore; we give its proof for the sake of
completeness.
\begin{theorem}\label{th:CD-arithm}
Assume that the distribution of $\beta$ is arithmetic, i.e.\ its
support $\varLambda$ is contained in the set $\lambda\ZZ$ with
maximal span $\lambda>0$. Then the general b.c.-solution of equation
\eqref{eq:CD1} is given by $y(x)=\g(x/\lambda)$, where $\g(\cdot)$
is any continuous periodic function of period~$1$.
\end{theorem}

\proof We start by showing that the smallest additive subgroup
$\calG\subset\lambda\ZZ$ generated by $\varLambda$
coincides with $\lambda\ZZ$.
Indeed, for $n\in\NN$ let $d_n\in\lambda\NN$ be the greatest common
divisor of the (finite) set
$\varLambda_n:=\{s\in\varLambda\colon\allowbreak |s|\le \lambda
n\}\subset \calG$. By B\'ezout's identity (see e.g.\
\cite[\S1.2]{JJ}) we have $d_n=\sum_{s_i\in\varLambda_n}\mynn m_i
s_i$ \strut{}with some integers $m_i$, and it follows that $d_n\in
\calG$ ($n\in\NN$). Since the sequence \strut{}$d_n/\lambda\in\NN$
is non-increasing, there exists the limit
$k^*\!:=\lim_{n\to\infty}d_n/\lambda=d_{n^*}/\lambda\in\NN$ and so
$\varLambda\subset k^* \lambda\ZZ$. But $\lambda>0$ is the maximal
span, hence $k^*=1$ and thus $\lambda=d_{n^*}\!\in \calG$, which
implies $\calG=\lambda\ZZ$, as claimed.

Now, it is easy to see that equation \eqref{eq:CD1} splits into
separate discrete equations $\tilde{z}(k)=\EE\{\tilde{z}(k-
\beta/\lambda)\}$ on every coset $x_0+\lambda\ZZ$
\,($x_0\in[0,\lambda)$), where $\tilde{z}(k):=z(x_0+k\lambda)$
($k\in\ZZ$). The discrete Choquet--Deny theorem shows that
$\tilde{z}(k)$ is constant on $\ZZ$; in other words, any bounded
solution of \eqref{eq:arch>0'} on $\RR$ is $\lambda$-periodic, and
the claim of the theorem easily follows.
\endproof

\begin{remark}
It is evident that any function $y(x)=\g(x/\lambda)$ satisfies
equation \eqref{eq:CD1}; the main point of
Theorem~\ref{th:CD-arithm} is that \emph{there are no other}
b.c.-solutions.
\end{remark}

\begin{remark}
Laczkovich
\cite{Laczkovich}
\strut{}gives a full characterization of
non-negative measurable solutions of the equation
$y(x)=\sum_{i=1}^\ell C_i\mypp y(x-b_i)$ (with arbitrary
coefficients $C_1,\dots,C_\ell>0$) in terms of the real roots of the
characteristic equation $\sum_{i=1}^\ell C_i \mypp \re^{-b_is}=1$.
\end{remark}

\begin{remark}
The original proof by Choquet \& Deny \cite{CD} (as well as many
subsequent proofs and extensions)
is based on a reduction to a uniformly continuous solution (cf.\
Remark~\ref{rm:UC}) and on establishing that the latter must reach
its maximum at a finite point $x_0\in\RR$. In view of the martingale
techniques used in the present paper for a general AE, it is of
interest to point out an elegant martingale proof found by Sz\'ekely
\& Zeng \cite{Szekely} (cf.\ Rao \& Shanbhag \cite[Ch.\;3]{Rao}).
\end{remark}

\subsubsection{Case $|\alpha|\equiv1$.} We prove here
an extension of the Choquet--Deny theorem for $\alpha$ taking the
values $\pm1$; to the best of our knowledge, such a result has not
yet been mentioned in the literature.
\begin{theorem}\label{th:CD}
Suppose that $|\alpha|\equiv 1$ and\/ $\PP(\alpha=1)<1$. Let
$\beta^+$, $\beta^-$ have the distribution of\/ $\beta$ conditioned
on $\alpha=1$ and $\alpha=-1$, respectively \textup{(}in case
$\alpha\equiv -1$, set $\beta^{+}\!\equiv0$\myp\textup{)}.

\begin{itemize}
\item[\rm (a)] \,If\/ $\beta^+\mynn$ is non-arithmetic then every b.c.-solution of
equation \eqref{eq:arch} is constant.

\item[\rm (b)]  \,Let
$\beta^+\mynn$ be arithmetic with span $\lambda\ne0$.
\begin{itemize}
\item[\rm (b-i)]
If\/ the distribution of\/ $\beta^-\mynn$ is not supported on any
set $\lambda_0+\lambda\ZZ$ \,\textup{($\lambda_0\in\RR$)}, then
every b.c.-solution of equation \eqref{eq:arch} is constant.

\item[\rm (b-ii)]
Otherwise, the general b.c.-solution of equation \eqref{eq:arch} is
of the form $y(x)=\g(x/\lambda)$, where $\g(\cdot)$ is a continuous
$1$-periodic function symmetric about point $x_0=\tfrac12\myp
\lambda_0/\lambda$, i.e. $\g(x_0-x)=\g(x_0+x)$
\,\textup{(}$x\in\RR$\textup{)}.
\end{itemize}
\end{itemize}
\end{theorem}

\begin{remark} In part (b-ii), all functions
with the required symmetry property may be represented (though not
uniquely) as $\g(x)=\g_0(x-x_0)+\g_0(-x+x_0)$, where $\g_0(\cdot)$
is an arbitrary (continuous) $1$-periodic function. It is
straightforward
to check that so constructed functions $y(x)=\g(x/\lambda)$ satisfy
equation \eqref{eq:arch} (with $|\alpha|\equiv 1$), but part (b-ii)
asserts that \emph{all} b.c.-solutions are contained in this class.
\end{remark}

\proof[Proof\/ of\/ Theorem\/ \textup{\ref{th:CD}}] Consider the
stopping time $\tau_1:=\inf\{n\ge 1\colon A_n=1\}$. With
$p_1:=\PP(\alpha=1)<1$ and $q_1:=\PP(\alpha=-1)=1-p_1>0$, the
distribution of $\tau_1$ is given by\footnote{These formulas include
the case $p_1=0$ (under the convention $0^0:=1$), whereby $\tau_1=2$
a.s.}
\begin{equation}\label{eq:tau-dist}
\PP(\tau_1=1)=p_1,\qquad \PP(\tau_1=n)=q_1^2\myp p_1^{n-2}\quad (
n\ge2),
\end{equation}
hence $\tau_1<\infty$ a.s. Since $A_{\tau_1}=1$ (a.s.), from
\eqref{eq:Y} we have $X_{\tau_1}=x-D_{\tau_1}$, and by
Lemma~\ref{lm:stop1}
\begin{equation}\label{eq:u=} y(x)=\EE\{y(x-D_{\tau_1})\},\qquad
x\in\RR.
\end{equation}
Now, in view of the Choquet--Deny theorem, we need to investigate
whether the random shift $D_{\tau_1}$ is non-arithmetic, i.e.
$\PP(D_{\tau_1}\!\in \lambda\ZZ)<1$ for any $\lambda\in\RR$.

Let $(\beta^+_n)$ and $(\beta^-_n)$ be two sequences of i.i.d.\
random variables each, with the same distribution as $\beta^+$ and
$\beta^-\mynn$, respectively. Conditioning on $\tau_1$ and using
\eqref{eq:AY} and \eqref{eq:tau-dist}, we obtain
\begin{align}
\notag \PP(D_{\tau_1}\!\in \lambda\ZZ)
\notag &=p_1\PP(\beta^+_1\in \lambda\ZZ) +\sum_{n=3}^\infty
q_1^2\myp
p_1^{n-2}\PP(\beta^-_1-\beta^+_2-\dots-\beta^+_{n-1}-\beta^-_n\in
\lambda\ZZ)\\
\label{eq:total} &\le p_1 +q_1^2\sum_{n=2}^\infty p_1^{n-2}=1.
\end{align}
If $p_1>0$ then \eqref{eq:total} implies that $\PP(D_{\tau_1}\!\in
\lambda\ZZ)<1$ unless $\beta_1^+\in\lambda\ZZ$ (a.s.) and for all
$n\ge 2$
\begin{equation}\label{eq:arith2}
\beta^-_1-\beta^+_2-\dots-\beta^+_{n-1}-\beta^-_n\in \lambda\ZZ
\quad\text{(a.s.)}.
\end{equation}
Since all $\beta_i^+$ are i.i.d., the first of these inclusions
implies that $\beta^+_2+\dots+\beta^+_{n-1}\in \lambda\ZZ$ a.s.\ for
all $n\ge2$; hence, conditions \eqref{eq:arith2} are reduced to
$\beta_1^{-}\!-\beta_2^{-}\!\in \lambda\ZZ$ (a.s.). In turn, the
last condition is equivalent to $\beta^-\!\in \lambda_0+\lambda\ZZ$
for some $\lambda_0\in\RR$. Indeed, applying Lebesgue's
decomposition theorem
\cite[p.\,142]{Feller2} to each of the i.i.d.\ random variables
$\beta_1^{-}$ and $\beta_2^{-}$, it is evident that the continuous
part of their common distribution must vanish, so that this
distribution is purely discrete; furthermore, its (countable)
support $\{b_i\}$ satisfies the condition $b_i-b_j\in\lambda\ZZ$ for
all $i,j$, and the claim follows.

A similar argument is also valid for $p_1=0$, whereby $\beta^+\!=0$,
\,$\beta^-\!=\beta$, and \eqref{eq:total} simplifies to
\begin{equation*}
\PP(D_{\tau_1}\!\in \lambda\ZZ)= \PP(\beta_1-\beta_2\in \lambda\ZZ)
\le 1.
\end{equation*}
This completes the proof of parts (a) and ~(b-i).

Finally, we prove part (b-ii), whereby $\beta^+\!\in\lambda\ZZ$,
$\beta^-\!\in \lambda_0+\lambda\ZZ$ and $D_{\tau_1}\!\in \lambda\ZZ$
(a.s.). By Theorem \ref{th:CD-arithm}, any b.c.-solution of equation
\eqref{eq:u=} must be of the form $y(x)=\g(x/\lambda)$, with some
$1$-periodic function $\g(\cdot)$. Substituting this into the
original equation \eqref{eq:u=} (with $|\alpha|\equiv1$) we get
\begin{align}
\notag
\g(x/\lambda)&=p_1\EE\{\g((x-\beta^+)/\lambda)\}+q_1\EE\{\g((-x+\beta^-)/\lambda)\}\\[.2pc]
\label{eq:q1>0}&=p_1 \g(x/\lambda)+q_1 \g((-x+\lambda_0)/\lambda),
\end{align}
and since $q_1\ne0$ it follows that $\g(\cdot)$ satisfies the
functional equation $\g(x)=\g(-x+\lambda_0/\lambda)$ \,($x\in\RR$),
which is equivalent to the symmetry condition stated in the theorem.
\endproof

\begin{remark}
If $\alpha\equiv1$ then $\beta^+\!\equiv\beta$,
\,$\beta^-\!\equiv0$, and parts (a) and (b-ii) of Theorem
\ref{th:CD} formally recover the Choquet--Deny theorem. No extra
requirement on $\g(\cdot)$ arises from \eqref{eq:q1>0}, since
$q_1=0$.
\end{remark}

\begin{remark}
Note that the values $\alpha=1$ and $\alpha=-1$ (and the
corresponding conditional distributions of $\beta$ represented by
$\beta^+\mynn$ and $\beta^-\mynn$, respectively) feature in Theorem
\ref{th:CD} in a non-symmetric way: e.g.\ if $\beta^+\myn$ has a
non-arithmetic distribution then, according to part (a), there are
no b.c.-solutions except constants, irrespectively of $\beta^-\myn$;
however, if $\beta^+\mynn$ is arithmetic with span $\lambda\ne0$
whilst $\beta^-\mynn$ is non-arithmetic but supported on a set
$\lambda_0+\lambda\ZZ$ (i.e.\ with $\lambda_0\ne0$ incommensurable
with $\lambda$) then there exist non-trivial b.c.-solutions,
according to part~(b-ii).
\end{remark}

\begin{example}\label{ex:CD}
Theorem \ref{th:CD} with $\alpha\equiv -1$ is exemplified by the
equations
\begin{itemize}
\item[(a)]
$y(x)=\int_0^\infty y(t-x)\,\re^{-t}\,\rd{t}$ \,(equivalent to the
pantograph equation $y'(x)+y(x)=y(-x)$, cf.~\eqref{eq:pant-gen}),
which by part (b-i) has only constant b.c.-solutions;
\item[(b)] $y(x)=\tfrac13\mypp y(-x+1)+\tfrac23\mypp y(-x-1)$,
which has periodic solutions of the form $y(x)={\rm g}(x)+{\rm
g}(-x)$, in accordance with part~(b-ii).
\end{itemize}
\end{example}

\subsection{The resonance case}\label{sec:resonance}

\begin{definition}\label{def:resonance}
The random coefficients $\alpha,\myp\beta$ of the AE \eqref{eq:arch}
are said to be \emph{in resonance} if there is a non-random constant
$c\in\RR$ such that $\alpha\myp(c-\beta)= c$
\,\textup{(}a.s.\textup{)}.
\end{definition}

The special role of resonance is clear from the observation that if
$X_0=c$ then by recursion \eqref{eq:Xn} we have $X_n= c$ (a.s.)\ for
all $n\ge 0$. It turns out that a Liouville-type theorem is always
true in the resonance case. Recall that we assume $\alpha\ne0$ a.s.
\begin{theorem}\label{th:degenerate}
Let $\PP(|\alpha|\ne1)>0$, and suppose that $\alpha,\,\beta$ are in
resonance. Then any b.c.-solution of the AE\/ \eqref{eq:arch} is
constant.
\end{theorem}
\proof Let $c\in\RR$ be such that $\alpha\myp(c-\beta)\equiv c$.
Observe that the substitution $\tilde{y}(x):=y(x+c)$ eliminates the
random shift in equation~\eqref{eq:arch},
\begin{align*}
\tilde{y}(x)=y(x+c) =\EE\{\alpha
x+\alpha\myp(c-\beta))\}=\EE\{y(\alpha\myp
x+c)\}=\EE\{\tilde{y}(\alpha\myp x)\}.
\end{align*}
Thus, without loss of generality, we can consider the equation
$y(x)=\EE\{y(\alpha\myp x)\}$. Denote $\tau_+:=\allowbreak\inf\{n\ge
1\colon A_n>0\}$ and $p:=\PP(\alpha>0)$, then
\begin{align*}
\PP(\tau_+>n)&= (1-p)\, p^{n-1}\to0\qquad (n\to\infty),
\end{align*}
so that $\tau_+<\infty$ \,a.s. Hence, by virtue of Lemma
\ref{lm:stop1} the equation $y(x)=\EE\{y(\alpha\myp x)\}$ is reduced
to $y(x)=\EE\{y(\tilde{\alpha}\myp x)\}$ with
$\tilde{\alpha}:=A_{\tau_+}>0$ \,a.s.

We first check that $\PP(\tilde{\alpha}=1)<1$. Indeed, with
$p=\PP(\alpha>0)$ as above we have
\begin{align}
\PP(\tilde{\alpha}=1) &=p\PP(A_1=1)+ \sum_{n=2}^\infty
(1-p)^2p^{n-2}\PP(A_{n}=1) \label{eq:tilde-alpha=1} \le p
+(1-p)^2\sum_{n=2}^\infty p^{n-2}=1.
\end{align}
If the probability on the left-hand side of \eqref{eq:tilde-alpha=1}
equals $1$ and $p>0$, then we must have
$\PP(A_1=1)=\PP(\alpha=1)=1$, which contradicts the theorem's
hypothesis; similarly, if $p=0$ then the inequality
\eqref{eq:tilde-alpha=1} implies that $A_2=\alpha_1\alpha_2=1$
(a.s.), and since $\alpha_1$, $\alpha_2$ are i.i.d.\ the latter
equality is possible only if $|\alpha|=1$ a.s., which is again a
contradiction.

Now, the equation $y(x)=\EE\{y(\tilde\alpha\myp x)\}$ on $\RR$ (with
$\tilde\alpha>0$) splits into two separate equations, for $x\ge0$
and $x\le0$, linked by the continuity condition at zero. For
instance, consider the equation
\begin{equation}\label{eq:arch>0}
y(x)=\EE\{y(\tilde\alpha\myp x)\} \quad (x>0), \qquad
y(0)=\lim_{x\to0+}y(x).
\end{equation}
Our aim is to show that $y(x)\equiv\const$ for all $x\ge0$. By the
change of variables
\begin{equation}\label{eq:z(t)}
t=-\ln x\in\RR,\qquad z(t)=y(\re^{-t}),
\end{equation}
the initial value problem \eqref{eq:arch>0} is transformed into
\begin{equation}\label{eq:arch>0'}
z(t)=\EE\{z(t-\tilde{\beta}\myp )\}\quad \ (t\in\RR),\qquad
\lim_{t\to+\infty}z(t)=y(0),
\end{equation}
which is an archetypal equation \eqref{eq:arch} with the unit
rescaling coefficient and random shift $\tilde\beta:=\allowbreak
\ln\tilde\alpha$ \,(such that $\PP(\tilde{\beta} \ne0)>0$), subject
to an additional limiting condition at $+\infty$.

If $\tilde{\beta}$ is non-arithmetic then the Choquet--Deny theorem
readily implies that all b.c.-solutions of equation
\eqref{eq:arch>0'} are constant (even without the limit condition at
$+\infty$).
In the arithmetic case, by Theorem \ref{th:CD-arithm} any
b.c.-solution of \eqref{eq:arch>0'} is $\lambda$-periodic, but due
to the limit in \eqref{eq:arch>0'} it must be constant. Returning to
\eqref{eq:arch>0} via the substitution \eqref{eq:z(t)}, we conclude
that in all cases $y(x)\equiv\const$ for $x\ge0$. By symmetry, the
same is true for $x\le0$, and the proof is completed by invoking
continuity of $y(x)$ at $x=0$.
\endproof

\begin{assumption}\label{as1}
Henceforth, unless explicitly stated otherwise, we assume that
$$
\text{(i)} \,\PP(\alpha\ne0)=1;\quad \text{(ii)}
\,\PP(|\alpha|\ne1)>0;\quad \text{(iii)} \ \,\alpha,\beta \
\,\text{are not in resonance}.
$$
\end{assumption}

\section{The critical case}\label{sec:3}

\subsection{Liouville theorem subject to uniform continuity}\label{sec:3.1}

Recall the notation $K:=\EE\{\ln\myn|\alpha|\}$. The next theorem
deals with the case $K\le0$ (including the critical case $K=0$)
under an additional \textit{a priori} hypothesis of uniform
continuity of the solution; on the other hand, in contrast to
Theorem \ref{th:K<>0} no moment conditions are imposed on $\beta$.
Note that conditions (i) and (ii) of Assumption \ref{as1} are in
force, but (iii) is not needed.

\begin{theorem}\label{th:K=0} Assume that\/ \/ $\PP(|\alpha|\ne
1\}>0$ and\/ $\EE\{\ln\myn|\alpha|\}\le 0$. Let $y(x)$ be a bounded
solution of the AE \eqref{eq:arch} which is uniformly continuous on
$\RR$. Then $y(x)\equiv \const$.
\end{theorem}
\proof By uniform continuity, for any $\varepsilon>0$ there is
$\delta=\delta(\varepsilon)>0$ such that if $|x_1-x_2|<\delta$ then
$|y(x_1)-y(x_2)|<\varepsilon$ \,($x_1,\myp x_2\in\RR$). Furthermore,
for a given $x\in\RR$ choose $M=M(\delta,x)>0$ such that
$|x|\,\re^{-M}<\delta$, and define the stopping time
\begin{equation*}
\tau_{M}:=\inf\{n\ge1\colon |A_n|\le \re^{-M}\}.
\end{equation*}
Since $\EE\{\ln\myn|\alpha|\}\le 0$, it follows (see e.g.\
\cite[pp.\;395--397]{Feller2}) that
\begin{equation*}
\liminf_{n\to\infty} \ln\myn|A_n|=\liminf_{n\to\infty} \sum_{k=1}^n
\ln\myn|\alpha_k|=-\infty\qquad \text{(a.s.)},
\end{equation*}
implying that $\tau_M<\infty$ a.s. Hence, by Lemma \ref{lm:stop1} we
have $y(x)=\EE_x\{y(X_{\tau_M})\}$ or, more explicitly (using the
iteration formulas \eqref{eq:Y}, \eqref{eq:AY}),
\begin{equation}\label{eq:E_x}
y(x)=\EE\{y(x A_{\tau_M}-D_{\tau_M})\},\qquad x\in\RR.
\end{equation}
In particular, \eqref{eq:E_x} with $x=0$ gives
$y(0)=\EE\{y(-D_{\tau_M})\}$. Hence, from \eqref{eq:E_x} we obtain
\begin{equation}\label{eq:d<Ed}
|y(x)-y(0)|\le \EE\bigl|y(x\myp
A_{\tau_M}-D_{\tau_M})-y(-D_{\tau_M})\bigr|.
\end{equation}
But according to the definition of the stopping time $\tau_M$ and
the choice of $M$, we have
$$
\bigl|(x\myp
A_{\tau_M}-D_{\tau_M})-(-D_{\tau_M})\bigr|=|x|\cdot|A_{\tau_M}|\le
|x|\,\re^{-M}<\delta.
$$
Due to uniform continuity of $y(x)$ (see above), this implies
$$
\bigl|y(x A_{\tau_M}-D_{\tau_M})-y(-D_{\tau_M})\bigr|<\varepsilon,
$$
and from \eqref{eq:d<Ed} we readily obtain $|y(x)-y(0)|\le
\varepsilon$. Since $\varepsilon>0$ is arbitrary, it follows that
$y(x)\equiv y(0)$, which completes the proof.
\endproof

\subsection{Sufficient conditions for the uniform continuity of solutions}\label{sec:3.2}


Consider the general AE \eqref{eq:arch} with no restriction on the
value $K=\EE\{\ln\myn|\alpha|\}$. Simple sufficient conditions for
the uniform continuity of its solutions are based on the following
well-known fact from real analysis (see e.g.\ \cite[p.\:74,
Proposition~2.5]{SS}).
\begin{lemma}\label{lm:SS}
Suppose that $f\in L^1(\RR)$. Then
\begin{equation*}
\lim_{h\to0}\int_\RR \myn|f(t+h)-f(t)|\,\rd{t}=0.
\end{equation*}
\end{lemma}
For continuous functions with compact support the lemma holds by
dominated convergence;
the general case follows since such functions are dense in
$L^1(\RR)$ \cite[p.\:71, Theorem 2.4(iii)]{SS}.

\begin{theorem}\label{th:uniform}
Assume that for each $a\in\RR$ in the support of the random variable
$\alpha$ there exists the conditional density function
$f_\beta(t\myp|\myp a)=\PP(\beta\in
\rd{t}\mypp|\mypp\alpha=a)/\rd{t}$. Then any bounded solution of\/
\eqref{eq:arch} is uniformly continuous on~$\RR$.
\end{theorem}

\proof Let $\mu_\alpha(\rd{a}):=\PP(\alpha\in\rd{a})$. By Fubini's
theorem, equation \eqref{eq:eq0} becomes
\begin{align*}
y(x)&=\int_{\RR} \mu_\alpha(\rd{a}) \int_\RR y(a\myp(x-t))\,
f_\beta(t\mypp|\mypp a)\,\rd{u}\equiv\int_{\RR} \mu_\alpha(\rd{a})
\int_\RR y(-au)\, f_\beta(x+u\mypp|\mypp a)\,\rd{u},
\end{align*}
where we used the change of variables $u=t-x$. Hence, uniformly in
$x\in\RR$
\begin{align*}
|y(x+h)-y(x)|&\le \int_{\RR} \mu_\alpha(\rd{a}) \int_\RR
|y(-au)|\cdot \bigl|f_\beta(x+h+u\mypp|\mypp a)-f_\beta(x+u\mypp|\mypp a)\bigr|\,\rd{u}\\
&\le \|y\| \int_\RR\mu_\alpha(\rd{a}) \int_\RR
\bigl|f_\beta(t+h\mypp|\mypp a)-f_\beta(t\mypp|\mypp
a)\bigr|\mypp\rd{t}\to0\qquad (h\to0),
\end{align*}
due to the bound $\|y\|<\infty$ and also using Lemma~\ref{lm:SS}
(applied to $f_\beta(\cdot\myp|\mypp a)$ for each $a$)
and Lebesgue's dominated convergence theorem.
\endproof

In many cases, another sufficient condition is more suitable. Let
the random variable $\beta$ in the AE \eqref{eq:arch} be of the form
$\beta=\gamma+\xi$, where $\xi$ is independent of the random pair
$(\alpha,\gamma)$ and has the density function, $f_\xi(t)$. Set
\begin{equation}\label{eq:phi1}
\varphi(x):=\EE\{y(\alpha\myp(x-\gamma))\}\equiv\iint_{\RR^2}
y(a\myp(x-c))\,\mu_{\alpha,\gamma}(\rd{a},\rd{c}),
\end{equation}
where
$\mu_{\alpha,\gamma}(\rd{a},\rd{c}):=\PP(\alpha\in\rd{a},\gamma\in\rd{c})$.
For example, if the measure $\mu_{\alpha,\gamma}$ is discrete, with
atoms $(a_i,c_i)$ and respective masses $p_i$, then
$$
\varphi(x)=\sum\nolimits_{i} p_{i}\mypp y(a_i(x-c_i)),\qquad
x\in\RR.
$$

Observe that, by independence of $\xi$, we have
\begin{align*}
\EE\{y(\alpha\myp(x-\beta))\}&=\EE\{y(\alpha\myp(x-\xi-\gamma))\}\\
&=\int_\RR \EE\{y(\alpha\myp(x-t-\gamma))\}\mypp\myp
f_\xi(t)\,\rd{t}=\int_\RR \varphi(x-t)\,f_\xi(t)\,\rd{t},
\end{align*}
according to the definition \eqref{eq:phi1}. Hence, the AE
\eqref{eq:arch} becomes
\begin{equation}\label{eq:4-}
y(x)=\int_\RR \varphi(x-t)\,f_\xi(t)\,\rd{t},\qquad x\in\RR.
\end{equation}
\begin{theorem}\label{th:uniform2}
Any bounded solution of equation\/ \eqref{eq:4-} is uniformly
continuous on~$\RR$.
\end{theorem}

\proof
By the substitution $u=x-t$ equation \eqref{eq:4-} is rewritten as
\begin{equation*}
y(x)=\int_\RR \varphi(u)\,f_\xi(x-u)\,\rd{u},\qquad x\in\RR.
\end{equation*}
Hence, uniformly in $x\in\RR$
\begin{align*}
|y(x+h)-y(x)| &\le \|\varphi\|\int_\RR \bigl|f_\xi(x+h-u)-
f_\xi(x-u)\bigr|\mypp\rd{u}\\
&= \|\varphi\| \int_\RR
\bigl|f_\xi(t+h)-f_\xi(t)\bigr|\,\rd{t}\to0\qquad (h\to0),
\end{align*}
according to Lemma~\ref{lm:SS}, and the claim follows.
\endproof

\begin{remark}
In both Theorems \ref{th:uniform} and \ref{th:uniform2}, continuity
of solutions is not assumed \emph{a priori}.
\end{remark}


\subsection{Pantograph equation}\label{sec:A1}
In this subsection, we explain the link pointed out in
\S\ref{sec:pantograph} between the AE \eqref{eq:arch} and a class of
functional-differential `pantograph' equations. We begin with an
elementary proof for a first-order pantograph equation
(\S\myp\ref{sec:A1.1}), and then treat the general case
(\S\myp\ref{sec:A1.2}). In turn, this allows us to establish the
uniform continuity of solutions by virtue of Theorem
\ref{th:uniform2}, and hence to prove a Liouville theorem using
Theorem~\ref{th:K=0} (\S\myp\ref{sec:A1.3}).

\subsubsection{First-order pantograph equation.}\label{sec:A1.1}
Assume that the random variable $\beta$ in \eqref{eq:arch} is
independent of $\alpha$ and has the unit exponential distribution,
with the density function $\re^{-t}\bfOne_{(0,\infty)}(t)$. Then
equation \eqref{eq:arch} specializes to
\begin{equation}\label{eq:4}
y(x)=\int_0^\infty \varphi_0(x-t)\,\re^{-t}\,\rd{t},\qquad x\in\RR,
\end{equation}
where
\begin{equation*}
\varphi_0(x):=\EE\{y(\alpha\myp x)\}\equiv\int_\RR y(a
x)\,\mu_\alpha(\rd{a}),\qquad
\mu_\alpha(\rd{a}):=\PP(\alpha\in\rd{a}).
\end{equation*}

\begin{proposition}\label{pr:pantograph}
Every b.c.-solution of\/ \eqref{eq:4} satisfies the pantograph
equation
\begin{equation}\label{eq:5}
y'(x)+y(x)=\varphi_0(x),\qquad x\in\RR.
\end{equation}
Conversely, any bounded solution of\/ \eqref{eq:5} satisfies
equation~\eqref{eq:4}.
\end{proposition}

\proof The substitution $u=x-t$ transforms equation \eqref{eq:4}
into $y(x)=\re^{-x}\!\int_{-\infty}^x
\varphi_0(u)\,\re^{u}\mypp\rd{u}$, and it is now evident that the
right-hand side is continuous and, moreover, \strut{}differentiable
in ${x\in\RR}$. Hence, by the Newton--Leibniz theorem
we readily obtain equation~\eqref{eq:5}.

Conversely, let $y(x)$ be a bounded solution of~\eqref{eq:5}. Then
by variation of constants
\begin{equation}\label{eq:6a}
y(x)= y(0)\,\re^{-x}+\int_0^x \varphi_0(u)\,\re^{\myp u-x}\,\rd{u}.
\end{equation}
Since $y(x)$ is bounded, we have $y(x)\mypp \re^{x}\to0$ as
$x\to-\infty$, and it follows from \eqref{eq:6a} that
$y(0)=\int_{-\infty}^0 \varphi_0(u)\,\re^{u}\,\rd{u}$. Substituting
this back into \eqref{eq:6a} and combining the integrals, we obtain
\begin{equation*}
y(x)=\int_{-\infty}^x
\varphi_0(u)\,\re^{u-x}\,\rd{u}=\int_0^{\infty}
\varphi_0(x-t)\,\re^{-t}\,\rd{t},
\end{equation*}
which is exactly equation~\eqref{eq:4}.
\endproof

\subsubsection{Higher-order pantograph equations.}\label{sec:A1.2}

The correspondence demonstrated in \S\myp\ref{sec:A1.1} can be
extended to more general equations, including higher orders. Like in
Theorem \ref{th:uniform2}, suppose that $\beta=\gamma+\xi$, where
$\xi$ is independent of $(\alpha,\gamma)$ and has density
$f_\xi(t)$. Following \cite{Der}, fix $r\in\NN$ and (real) constants
$\kappa_1,\dots,\kappa_r\ne0$ (some or all of which may coincide),
and let $f_\xi(t)$ be given by
\begin{equation*}
f_\xi(t)=g_{1}\star\dots\star g_{r}(t),\qquad g_{j}(t):=|\kappa_j|\,
g_0(\kappa_j\myp t)\quad(j=1,\dots,r),
\end{equation*}
where $\star$ denotes convolution and
$g_0(t):=\re^{-t}\bfOne_{(0,\infty)}(t)$. On the other hand,
consider the pantograph equation of order~$r$
\begin{equation}\label{eq:Q}
\prod_{j=1}^r
\left(1+\frac{\calD}{\kappa_j}\right)y(x)=\varphi(x),\qquad
\calD:=\frac{\rd}{\rd{x}},
\end{equation}
where (see \eqref{eq:phi1})
\begin{equation}\label{eq:phi2}
\varphi(x)=\EE\{y(\alpha\myp(x-\gamma))\},\qquad x\in\RR.
\end{equation}
Finally, recall from \eqref{eq:4-} that the AE \eqref{eq:arch} is
equivalently rewritten as
\begin{equation}\label{eq:4+}
y(x)=\int_\RR \varphi(x-t)\,f_\xi(t)\,\rd{t},\qquad x\in\RR.
\end{equation}

\begin{proposition}\label{pr:pantograph-n}
Every b.c.-solution of \eqref{eq:4+} satisfies the pantograph
equation~\eqref{eq:Q}. Conversely, any bounded solution of
\eqref{eq:Q} satisfies equation~\eqref{eq:4+}.
\end{proposition}

\proof Let $y(x)$ be a b.c.-solution of \eqref{eq:4+}. Then the
function $\varphi(x)$ defined in \eqref{eq:phi2} is bounded,
$\|\varphi\|\le \|y\|<\infty$; furthermore, it is continuous by
virtue of dominated convergence.
Therefore, one can
apply the convolution inversion formula (see \cite[
p.\;28, Theorem~7.1]{HW}), which readily yields that $y(x)$
satisfies equation~\eqref{eq:Q}.

If $y(x)$ is a bounded (continuous) solution of \eqref{eq:Q}, then
$\varphi(x)$ is again b.c.\ and it follows that $y(x)$ satisfies
equation \eqref{eq:4+} thanks to \cite[p.\;34, Theorem~9.3]{HW}.
\endproof

\begin{example}\label{ex:pant}
The following are three simple examples illustrating
Proposition~\ref{pr:pantograph-n}.
\begin{itemize}
\item[(a)]
Take $r=1$ and $\kappa_1=1$, then equation \eqref{eq:Q} is reduced
to $y'(x)+y(x)=\varphi(x)$ (cf.~\eqref{eq:5}). Here, $\xi$ has the
unit exponential density, $f_\xi(t)=g_{0}(t)$. Thus, Proposition
\ref{pr:pantograph-n} extends Proposition \ref{pr:pantograph} to
pantograph equations with the right-hand side $\varphi(x)$ given
by~\eqref{eq:phi2}.

\item[(b)]
For $r=2$ and $\kappa_1=1$, $\kappa_2=-1$, equation \eqref{eq:Q}
specializes to $-y''(x)+y(x)=\varphi(x)$. This corresponds to a
symmetric (two-sided) exponential density $f_\xi(t)=g_{0}(t)\star
g_{0}(-t)=\frac12\, \re^{-|t|}$ ($t\in\RR$).

\item[(c)]
Likewise, with $r=2$ and $\kappa_1=\kappa_2=1$ we have $y''(x)+2\myp
y'(x)+y(x)=\varphi(x)$. Here, $f_\xi(t)=g_{0}\star
g_{0}(t)=t\,\re^{-t}$ ($t>0$).
\end{itemize}
\end{example}


\subsubsection{Liouville theorem for the pantograph
equation.}\label{sec:A1.3}

By virtue of Proposition \ref{pr:pantograph-n},
Theorem~\ref{th:uniform2} implies the following

\begin{proposition}\label{pr:UC} Any bounded
solution of the pantograph equation \eqref{eq:Q} is uniformly
continuous.
\end{proposition}

The next result settles a Liouville theorem for a general class of
pantograph equations in the critical and subcritical cases (cf.\
Theorem~\ref{th:pantograph1} stated in \S\ref{sec:intro2}).
\begin{theorem}\label{th:pantograph}
If $K=\EE\{\ln\myn|\alpha|\}\le 0$ then any bounded solution of
equation \eqref{eq:Q} is constant.
\end{theorem}

\proof By Proposition \ref{pr:UC}, bounded solutions of equation
\eqref{eq:pant-gen-alpha} are uniformly continuous; thus, if
$\PP(|\alpha|=1)<1$ then the claim readily follows by
Theorem~\ref{th:K=0}. If $|\alpha|=1$ a.s.\ then we can apply the
generalized Choquet--Deny theorem (Theorem \ref{th:CD}), noting that
either $\beta^+\mynn$ or $\beta^-\mynn$ (i.e.\ $\beta$ conditioned
$\alpha=1$ or $\alpha=-1$, respectively) must have a continuous
distribution, because $\beta=\gamma+\xi$ and $\xi$ is exponentially
distributed independently of $(\alpha, \gamma)$. Thus, the theorem
is proved.
\endproof



\subsection{$q$-Difference equations with shifts}\label{sec:3.3}

For the purely functional equation (cf.~\eqref{eq:pure1})
\begin{equation}\label{eq:16}
y(x)=\sum_{i=1}^\ell p_{i} \mypp y(a_i(x-b_i))\ \ \quad
(x\in\RR),\qquad p_i>0,\quad \sum_{i=1}^\ell p_i=1,
\end{equation}
where the above criteria of uniform continuity (Theorems
\ref{th:uniform},~\ref{th:uniform2}) are not usable, no general
results of the Liouville type are currently available in the
critical case, except for equations with $\alpha\equiv 1$ treated by
the Choquet--Deny theorem (see \S\S\myp\ref{sec:1.2.1},
\ref{sec:CD}) and its generalization to the case $|\alpha|\equiv 1$
(Theorem~\ref{th:CD}), and also equations without shifts (i.e.
$\beta\equiv0$) covered by Theorem~\ref{th:degenerate}.



In this subsection, we consider an important subclass of functional
equations \eqref{eq:16}, for which a Liouville theorem (in the
critical case) can be proved without the \emph{a priori} hypothesis
of uniform continuity. Namely, assume that the coefficients $a_i>0$
in \eqref{eq:16} are multiplicatively commensurable,
i.e. $a_i=q^{m_i}$ with some $q>1$ and $m_i\in\ZZ$
\,$(i=1,\dots,\ell)$. The resulting equations $y(x)=\sum_{i=1}^\ell
p_{i} \mypp y(q^{m_i}(x-b_i))$ are known as \textit{$q$-difference
equations}; the general theory of such equations (albeit without
shifts $b_i$) was developed by Birkhoff \cite{B} and Adams~\cite{A}.

To avoid trivialities, we assume in \eqref{eq:16} that $\ell\ge2$
and $(a_i,b_i)\ne(1,0)$ for all $i=1,\allowbreak \dots,\ell$. We
also exclude the case $a_1=\dots=a_\ell=1$, which is covered by the
Choquet--Deny theorem (see~\eqref{eq:CD}). The theorem below handles
the critical case, $K={\sum_{i=1}^\ell p_i\ln a_i = 0}$.

\begin{theorem}\label{th:3}
Assume that $\sum_{i=1}^\ell p_i\myp m_i =0$. Then any b.c.-solution
of equation \eqref{eq:16} is constant.
\end{theorem}
\proof Set
\begin{equation}\label{eq:gamma}
\rho_i:=\frac{b_i}{1-a_i^{-1}},\qquad i=1,\dots,\ell.
\end{equation}
If $a_i=1$ (but $b_i\ne0$) then \eqref{eq:gamma} is understood as
$\rho_i:=\infty$. Note that if $\rho_1=\dots=\rho_\ell=c\in\RR$
then equations \eqref{eq:gamma} are combined as
$\alpha\myp(c-\beta)\equiv c\myp$; that is to say, $\alpha$ and
$\beta$ are in resonance (see Definition \ref{def:resonance}), and
the desired result readily follows by Theorem~\ref{th:degenerate}.

Assuming now that not all $\rho_i$ are the same, let us follow a
similar martingale strategy as in the proof of Theorems
\ref{th:alpha=0}, \ref{th:CD} and \ref{th:K=0}, but based on the
stopping time $\tau_0:=\inf\{n\ge1\colon S_n=0\}$ of the random walk
$S_n:=\sum_{j=1}^n \log_q\alpha_j$ \,($S_0=0$). Since $\EE\{\log_q
\alpha\}=\sum_{i=1}^\ell p_i\myp m_i=0$, the random walk $S_n$ is
recurrent
(see e.g.\ \cite[p.\;33, Theorem~T1]{Spitzer}), hence
$\tau_0<\infty$ a.s. Thus, by Lemma \ref{lm:stop1}
\begin{equation}\label{eq:16.5}
\EE\{y(x-D_{\tau_0})\}=y(x),\qquad x\in\RR,
\end{equation}
where the sequence $(D_n)$ is defined in~\eqref{eq:AY}. By the
Choquet--Deny theorem (see~\S\myp\ref{sec:CD}), every b.c.-solution
of \eqref{eq:16.5} is constant on $\RR$ provided that \textit{the
distribution of the shift\/ $D_{\tau_0}$ is non-arithmetic}. In the
rest of the proof, our aim is to verify the last condition.

Due to recurrence of the random walk $(S_n)$, there are integers
$k_1,\dots,k_\ell\ge 0$ such that
\begin{equation}\label{eq:=0}
k_1m_1+\dots+k_\ell\mypp m_\ell=0,\qquad
k_1|m_1\mynn|+\dots+k_\ell\mypp |m_\ell|>0.
\end{equation}
This corresponds to paths $\{S_j\not\equiv0, \,0\le j\le
k^*\!:=k_1+\dots+k_\ell\}$ with exactly $k_i$ steps
$\log_q\alpha=m_i$ \,($i=1,\dots,\ell$), so that $S_{k^*}=0$. Let us
choose a sequence of these steps
so that $S_j>0$ for all $j=1,\dots,\myp k^*\!-1$, which would ensure
that $k^*$ is a possible value of the return time $\tau_0$ occurring
with probability $\PP(\tau_0=k^*)=p_1^{k_1}\!\cdots
p_\ell^{k_\ell}>0$. To this end, split the indices $i=1,\dots,\ell$
into disjoint groups with the same value $\rho_i$ in each group. Due
to the balance condition in \eqref{eq:=0}, the integers $m_i$ cannot
all have the same sign, and recalling that not all $\rho_i$
coincide, it is easy to see that one can find two indices $i^*$ and
$j^*$ such that
$$
\rho_{i^*}\ne \rho_{j^*}\mynn,\ \ \quad m_{i^*}\!>0,\ \ \quad
m_{j^*}\!<0,
$$
which also implies that $\rho_{i^*}$, $\rho_{j^*}$ are finite.
Hence, by a suitable relabelling of $a_1,\dots,a_\ell$ (so that
$i^*$ and $j^*$ become $1$ and $\ell$, respectively), we can assume
without loss of generality that $\rho_1\ne\rho_\ell$ and
$$
m_1>0,\ \ \quad m_2,\dots,m_r\ge0,\ \ \quad
m_{r+1},\dots,m_{\ell-1}\le 0,\ \ \quad m_\ell<0.
$$
In particular, it follows that
\begin{equation}\label{eq:sj}
s_i:=k_1 m_1+\dots+k_i\myp m_i>0,\qquad i=1,\dots,\myp \ell-1.
\end{equation}

Now, recalling that $a_1^{k_1}\!\cdots\myp
a_\ell^{k_\ell}=q^{k_1m_1+\dots+k_\ell\myp m_\ell}=q^0=1$, we obtain
\begin{equation}\label{eq:theta}
\theta_n:= \sum_{i=1}^\ell
\rho_i\myp(1-a_i^{-nk_i})\prod_{j=i}^{\ell}a_{j}^{nk_{j}}
=\sum_{i=1}^\ell
\rho_i\myp(1-a_i^{-nk_i})\prod_{j=1}^{i-1}a_{j}^{-nk_{j}},\qquad
n\in\NN.
\end{equation}
Equation \eqref{eq:theta} means that $\theta_n$ belongs to the
support of the random variable $D_{\tau_0}$.
Expanding and rearranging \eqref{eq:theta}, and using conditions
\eqref{eq:sj}, we have
\begin{align*}
\theta_n&=\sum_{i=1}^{\ell} \rho_i
\prod_{j=1}^{i-1}a_{j}^{-nk_{j}}-\sum_{i=1}^{\ell} \rho_i
\prod_{j=1}^{i}a_{j}^{-nk_{j}}\\
&=\rho_1+\sum_{i=1}^{\ell-1} (\rho_{i+1}-\rho_i)
\prod_{j=1}^{i}a_j^{-nk_{j}}-\rho_\ell
\prod_{j=1}^{\ell}a_j^{-nk_{j}}\\
&=\rho_1+\sum_{i=1}^{\ell-1} (\rho_{i+1}-\rho_i)\,q^{-n
s_i}-\rho_\ell\to \rho_1-\rho_\ell\ne0\qquad (n\to\infty).
\end{align*}
Hence, for any $\varepsilon>0$ there are $n,n'\in\NN$ such that
$0\ne |\theta_n-\theta_{n'}|<\varepsilon$, which implies that the
distribution of $D_{\tau_0}$ is non-arithmetic, as required.
\endproof

\appendix

\section{Appendix}\label{sec:A}

\subsection{Direct proof of Lemma~\ref{lm:stop1}}\label{sec:lemma}
The idea of the proof of the identity \eqref{eq:stop_x} is to
propagate equation \eqref{eq:Markov,n=1} according to the (random)
value of stopping time $\tau$, resulting in the relation
\begin{equation}\label{eq:tau>n*}
y(x)=\sum_{i=1}^n \EE_x\bigl\{y(X_i)\myp
I_{\{\tau=i\}}\bigr\}+\EE_x\bigl\{y(X_{n})\myp
I_{\{\tau>n\}}\bigr\},\qquad n\in\NN,
\end{equation}
where $I_A$ is the indicator of event $A$ (i.e.\ $I_A=1$ if $A$
occurs and $I_A=0$ otherwise).

We prove formula \eqref{eq:tau>n*} by mathematical induction. For
$n=1$ it is reduced to \eqref{eq:Markov,n=1}.
Suppose now that \eqref{eq:tau>n*} holds for some $n\in\NN$. For
shorthand, denote $\calF_n:=\sigma\{(\alpha_k,\beta_k),\,k\le n\}$,
i.e.\ the smallest $\sigma$-algebra containing all events
$\{\alpha_k\le a,\beta_k\le b\}$ ($a,b\in\RR$, $k=1,\dots,n$). The
assumption of the lemma that $\tau$ is a stopping time (relative to
$(\calF_n)$) means that $\{\tau\le n\}\in \calF_n$ ($n\in\NN$). By
the total expectation formula \cite[\S\myp{}9.7(a)]{Williams} we can
write
\begin{align}
\notag
\EE_x\{y(X_{n+1})\myp I_{\{\tau>n\}}\}&=\EE_x\bigl\{\EE_x\bigl[y(X_{n+1})\myp I_{\{\tau>n\}}|\mypp\calF_{n}\bigr]\bigr\}\\
&=\EE_x\bigl\{I_{\{\tau>n\}}\EE_x\bigl[y(X_{n+1})|\mypp\calF_{n}\bigr]\bigr\}
=\EE_x\bigl\{y(X_n)\myp I_{\{\tau>n\}}\bigr\}, \label{eq:tau>n**}
\end{align} thanks to the martingale
property~\eqref{eq:y(X)mart}; note that
$I_{\{\tau>n\}}=1-I_{\{\tau\le n\}}$ is $\calF_{n}$-measurable,
so can be taken outside the conditional expectation
$\EE_x[\dots\myn|\mypp\calF_{n}]$ \cite[\S9.7(j)]{Williams}.
Now, substituting equality \eqref{eq:tau>n**} into \eqref{eq:tau>n*}
and noting that
\begin{equation*}
\EE_x\{y(X_{n+1})\myp I_{\{\tau>n\}}\}=\EE_x\bigl\{y(X_{n+1})\myp
I_{\{\tau=n+1\}}\bigr\}+\EE_x\bigl\{y(X_{n+1})\myp
I_{\{\tau>n+1\}}\bigr\},
\end{equation*}
we conclude that the induction step is complete, and hence
\eqref{eq:tau>n*} holds for all $n\in\NN$.

Finally, since $\|y\|:=\sup_{x\in\RR} |y(x)|<\infty$ and
$\tau<\infty$ a.s., we have
$$
\bigl|\EE_x\bigl\{y(X_{n})\myp I_{\{\tau>n\}}\bigr\}\bigr|\le
\|y\|\PP_x(\tau>n)\to 0,\qquad n\to\infty.
$$
Hence, passing to the limit in \eqref{eq:tau>n*} as $n\to\infty$
yields
$$
y(x)=\sum_{i=1}^\infty \EE_x\{y(X_i)\myp
I_{\{\tau=i\}}\}=\EE_x\{y(X_\tau)\myp
I_{\{\tau<\infty\}}\}=\EE_x\{y(X_\tau)\},
$$
and the identity \eqref{eq:stop_x} is proved.

\subsection{Harmonic functions and asymptotics of Markov chains}\label{sec:A2}
Here, we provide a brief summary of some facts pertaining to the
role of harmonic functions in the asymptotic classification of
general Markov chains. For a more systematic exposition, we refer
the reader to monographs by Revuz~\cite{Revuz} and Meyn \& Tweedie
\cite{Meyn}.

Consider a general (time-homogeneous) Markov chain $(X_n)$ with
state space $\RR$. Denote by $\PP_x$ the probability law of $(X_n)$
started from $X_0=x\in\RR$, and by $\EE_x$ the corresponding
expectation. Let $\RR^\infty$ be the space of real sequences
$(x_0,x_1,\dots)$ and \strut{}$\calB^\infty$ the smallest
$\sigma$-algebra containing all cylinder sets
$\bar{x}=\{x\in\RR^\infty\colon (x_0,\dots,x_{m})\in B_{m}\}$ with a
Borel base $B_{m}\subset \RR^{m+1}$. Event $A$ is said to be
\emph{invariant} if there exists $\bar{B}\in\calB^\infty$ such that
$A=\{(X_m,X_{m+1},\dots)\in \bar{B}\}$ for every $m\ge 0$. For
example, the event $\{\lim_{n\to\infty} X_n=\infty\}$ is clearly
invariant because it is not affected by the shifts $(X_n)\mapsto
(X_{n+m})$ ($m\in\NN$). The class of all invariant events is a
$\sigma$-algebra denoted $\calI$; a random variable $Y$ is called
\emph{invariant} if it is $\calI$-measurable.

Recall that a function $y(x)$ is called \emph{harmonic} if
$y(x)=\EE_x\{y(X_1)\}$ for all $x\in\RR$.
A fundamental result (see e.g.\ \cite[
p.\;56, Proposition 3.2]{Revuz} or \cite[p.\;425, Theorem
17.1.3]{Meyn}) is that there is a one-to-one correspondence between
bounded harmonic functions and (equivalence classes of) bounded
invariant random variables, expressed by the pair of relations
\begin{equation}\label{eq:xi}
Y=\lim_{n\to\infty} y(X_n) \quad \text{(a.s.)},\qquad
y(x)=\EE_x\{Y\}\quad(x\in\RR).
\end{equation}
This result highlights the significance of information about bounded
harmonic functions in the asymptotic characterization of Markov
chains, especially with regard to transience \emph{vs.}\ recurrence
\cite[\S2.3]{Revuz}. For a Borel set $B\subset\RR$, put
$N_B:=\sum_{n=0}^\infty \bfOne_B(X_n)$ (\,$=\text{}$the total number
of visits to~$B$). Since $\{N_B=\infty\}\in\calI$, by \eqref{eq:xi}
the function $h_B(x):=\PP_x(N_B=\infty)$ is (bounded) harmonic. We
say that $B$ is \emph{transient} if $h_B(x)\equiv 0$ and
\emph{recurrent} (or \emph{Harris recurrent} \cite[Ch.~9]{Meyn}) if
$h_B(x)\equiv 1$. The following criterion is valid \cite[p.\;58,
Proposition 3.8]{Revuz}: in order that \emph{every (Borel) set be
either transient or recurrent}, it is necessary and sufficient that
either of the two equivalent conditions hold: (i) $\calI$ is trivial
(up to a.s.-equivalence); (ii) all bounded harmonic functions are
constant.

This criterion can be illustrated by our result on the existence of
a non-constant (continuous) solution of the AE \eqref{eq:arch} in
the supercritical case $K>0$, described in Theorem \ref{th:K<>0}(b).
Namely, take $B=(b,\infty)$ ($b\in\RR$), then by inspection of the
proof (see details in \cite{BDM}) it is evident that the function
$h_B(x)$ coincides with the solution $F_\Upsilon(x)=\PP(\Upsilon\le
x)$, where the random series $\Upsilon$ is defined in the theorem;
thus, any such set $B$ is neither transient nor recurrent.

\subsubsection*{Acknowledgements.} We are grateful to John Ockendon
and Anatoly Vershik for stimulating discussions, and to Charles
Taylor for helpful remarks. Hospitality and support during various
research visits provided by the Center for Advanced Studies in
Mathematics (Ben-Gurion University) and by the ZiF and SFB\,701
(Bielefeld University) are thankfully acknowledged. We also thank
the anonymous referees for the careful reading of the manuscript and
for the useful comments that have led to several improvements in the
presentation including the statements and proofs in \S\ref{sec:3.2}.

\subsubsection*{Funding statement.}
L.{}V.\,B.\ was partially supported by a Leverhulme Research
Fellowship; \,G.\,D.\ received funding from the Israel Science
Foundation (Grant 35/10).


\end{document}